\documentclass{amsart} 

\newtheorem{definition}{Definition}[section]

\newtheorem{remark}[definition]{Remark}
\newtheorem{example}[definition]{Example}
\newenvironment{ack}{\noindent{\bf Acknowledgements}.}{}

\newtheorem{lemma}[definition]{Lemma}
\newtheorem{proposition}[definition]{Proposition}
\newtheorem{theorem}[definition]{Theorem}
\newtheorem{corollary}[definition]{Corollary}

\def\P{{\mathbb{P}}}
\def\H{{\mathrm{H}}}
\def\K{{\mathbb{K}}}
\def\kK{{\mathcal{K}}}
\def\N{{\mathbb{N}}}
\def\CC{{\mathbb{C}}}

\def\I{{\mathcal{I}}}

\def\t{{\underline{t}}}
\def\s{{\underline{s}}}
\def\X{{\underline{X}}}
\def\u{{\underline{u}}}
\def\C{{\mathrm{C}}}
\def\S{{\mathrm{S}}}
\def\M{{\underline{M}}}
\def\ttheta{{\underline{\theta}}}
\def\tmu{{\underline{\mu}}}
\def\bdeg{{\mbox{bideg}}}
\begin{document}
\title[Generators of the Rees Algebra of monoid parametrizations]{Minimal generators of the defining ideal of the Rees Algebra associated to monoid 
parametrizations}

\author{Teresa Cortadellas Ben\'itez}
\address{Universitat de Barcelona, Departament d'{\`A}lgebra i Geometria.
Gran Via 585, E-08007 Barcelona, Spain}
\email{terecortadellas@ub.edu}

\author{Carlos D'Andrea}
\address{Universitat de Barcelona, Departament d'{\`A}lgebra i Geometria.
Gran Via 585, E-08007 Barcelona, Spain}
\email{cdandrea@ub.edu}
\urladdr{http://atlas.mat.ub.es/personals/dandrea}
\thanks{Both authors are supported by the Research Project MTM2007--67493 from the
Ministerio de Ciencia e Innovaci\'on, Spain}


\begin{abstract}
We describe a minimal set of generators of the defining ideal of the Rees Algebra associated to a proper
parametrization of any monoid hypersurface. In the case of plane curves, we recover a known description for
rational parametrizations having a syzygy of minimal degree ($\mu=1$). We also show that our approach can be applied
to parametrizations of rational surfaces having a Hilbert-Burch resolution with $\mu_1=\mu_2=1$. 
\end{abstract}
\maketitle

\section{Introduction}
In recent years a lot of attention has been given to the mathematical study of the so called ``method of implicitization
by  moving curves and surfaces'', which was introduced in the geometric modeling community by Sederberg and Chen 
in \cite{SC95} and subsequently studied and extended by several authors, see for instance 
\cite{SGD97,cox01,BCD03,CGZ00,CCL05} and the references therein.
\par A brief idea of this method is as follows. Let $n$ be a positive number and $\K$ a field. Consider the following
rational map
$$\begin{array}{ccc}
\P^{n-1}&\stackrel{\phi}{\dasharrow}&\P^{n}\\
(t_1:\ldots :t_n)&\mapsto&\big(u_1(\t):\ldots :u_{n+1}(\t)\big),   
  \end{array}
$$
where $\P^k$ denotes the $k$-dimensional projective space over $\K$, $\t=(t_1,\ldots,t_n)$ is a set of variables,
and $u_i(\t),\,i=1,\ldots, n+1,$ are homogeneous polynomials of degree $d>0$ in $\K[\t]$. The implicitization problem consists in finding
equations for the image of $\phi$ in $\P^{n}$. From a theoretical point of view, the solution of this problem is
quite easy as we are looking for generators of the elimination ideal
$$\langle X_1-u_1(\t),\,\ldots\,X_{n+1}-u_{n+1}(\t)\rangle\cap\K[\X],
$$
with $\X=(X_1,\ldots,X_{n+1}),$ and this task can be done with classical tools from elimination theory like 
 Gr\"obner bases or resultants (see \cite{CLO07} for more on this subject). However, in practice these
computations are very expensive. The method of moving surfaces exploits the idea of looking for
polynomials  $F(\t,\X)\in\K[\t,\X]$  such that $F(\t,\u(\t))=0$, having lower degree than $d=\deg(X_i-u_i(\t))$ . 
We can then recover the implicit equations by eliminating $\t$ from an appropriate family of these $F(\t,\X)$'s.
\par
In \cite{cox08}, this method was connected with the computation of generators of the defining ideal of
the  Rees Algebra  associated to the map $\phi$. Algebraically, the problem gets translated into the computation
of generators of the kernel of
the following $\K[\t]$-morphism of algebras:
$$\begin{array}{ccccc}
{\mathfrak h}:&\K[\t,\X]&\to&\K[\t][Z]&\\
&X_i&\mapsto& u_i\,Z\ &\quad i=1,\ldots, n+1.   
  \end{array}
$$
Here, $Z$ is a new variable. The Rees Algebra associated to $\phi$ is the image of ${\mathfrak h}$, and the
polynomials $F(\t,\X)$'s can be chosen among generators of $\mbox{ker}({\mathfrak h}).$  
This dictionary between moving surfaces and Rees Algebras had already been exploited by 
\cite{BJ03,BC05,BCJ09}, and proved to be fruitful, as shown by the subsequent papers \cite{bus09,CHW08,HSV08,HSV09, HW09, KPU09}.
\par
In this paper we give an elementary approach to the computation of minimal generators of $\mbox{ker}(\mathfrak{h})$ 
in the following cases:
\begin{itemize}
 \item $\phi$ a monoid parametrization,
\item $\phi$ a parametrization given by a sequence $u_1(\t),\,u_2(\t),\,u_3(\t),\,u_4(\t)$ 
($n=3$) which is a local complete intersection in $\P^2$ satisfying also that the ideal generated by the $u_i(\t)$'s 
is saturated and its syzygy module has two elements of degree $1$. In the language
of $\mu$-bases (see \cite{cox01,BCD03}), this means that $\tmu=(1,1,d-2)$.
\end{itemize}
A monoid parametrization is a proper map that parameterizes a so-called \textit{monoid hypersurface}, which are
rational varieties having one singularity of maximal rank (see Section 
\ref{s2} for the proper definition and references). For plane curves, being a monoid curve  is equivalent to the fact 
that there is a homogeneous element of bidegree $(1,1)$ in $\mbox{ker}(\mathfrak{h})$. In the language of $\mu$-bases (see 
\cite{CGZ00,cox01,bus09})), being a monoid curve means $\mu=1$. 
This was the case of study in \cite{cox08,CHW08,bus09,KPU09}. We recover their results in
Section \ref{s1} with elementary arguments (no need to introduce Sylvester Forms or Local Cohomology), 
and also show that the condition of having $\mu=1$ already implies that the map 
$\phi$ is generically injective (see Corollary \ref{cor}), and we can actually produce very easily an inverse of the
parametrization via the generator of the lowest degree part of $\mbox{ker}(\mathfrak{h}).$
\par
This is the key idea behind all the cases treated in this note: 
if we have enough elements of minimal degree in $\mbox{ker}({\mathfrak h})$, then we can read the inverse of 
$\phi$ from the part of lowest degree
in the defining ideal of the Rees Algebra associated to $\phi.$ Another key fact is that the ideal generated by the forms of lowest
degree is prime in $\K[\t,\X].$
We believe that our results can be generalized to a more general situation provided
that these two properties hold.
\par
We would like to point out that in \cite{HW09}, similar results are obtained for the case $n=3,\,\tmu=(1,1,d-2)$ by applying fine tools from local cohomology theory.
\par The paper is organized as follows: in Section \ref{s1} we deal with the case of parametrizations of plane
curves with $\mu=1$. 
Almost no new results are presented here, but the proofs are rather elementary and may help as a warm up towards 
the other cases. In Section \ref{s2} we deal with proper parametrizations $\phi$ of a general monoid
hypersurface and construct a minimal set of generators of $\mbox{ker}(\mathfrak{h})$ in this case. 
We show in Section \ref{s3} that our approach also works in the case of surfaces, provided that the ideal 
generated by the parametrization is saturated, local
complete intersection with two syzygies linearly  of degree $1$ (i.e. $\tmu=(1,1,d-2)$). 
Some technical results of Commutative Algebra will be used at the end of this section in order to prove the main result. 

\smallskip
\begin{ack}
We are grateful to David Cox for calling our attention to the paper \cite{KPU09}, and to J. William Hoffman for doing the same with \cite{HW09}.
We would also like to thank  Laurent Bus\'e for patiently explaining to us his recent work related
to this subject, and to the referees for their comments and suggestions in order to improve the presentation of this manuscript. All our experiments with computers were done with the aid of {\tt Maple} and {\tt Macaulay 2}.
\end{ack}

\section{Monoid Parametrizations of Plane Curves}\label{s1}

Let $\K$ be an algebraically closed field and 
${\mathcal I}:=\langle u_1(t_1,t_2),u_2(t_1,t_2),u_3(t_1,t_2)\rangle$ a
homogeneous ideal in $\K[t_1,t_2]$. Assume that $u_1,u_2,u_3$ are nonzero homogeneous polynomials of degree $d$ 
without common  factors. We are interested in computing a minimal set of generators of the kernel of the graded
morphism of $\K[t_1,t_2]$-algebras
$$\begin{array}{cccc}
   {\mathfrak h}:&\K[t_1,t_2,X_1,X_2,X_3]&\to&\mbox{Rees}({\mathcal I})\\
&X_i&\mapsto&u_i\,Z
  \end{array}
$$
for $i=1,2,3$. Here, $Z$ is a new variable, and $\mbox{Rees}(\I)=\K[t_1,t_2][\I\,Z]$ is the
Rees Algebra associated to $\I.$ Let $\kK\subset R[X_1,X_2,X_3]$ be the kernel of $\mathfrak{h}$. 
\par
Set $\t:=(t_1,t_2),\,\u:=(u_1,u_2,u_3)$ and $\X:=(X_1,X_2,X_3)$.  We have then that $\kK$ is the bi-graded ideal 
(with grading given by total degrees in $\t$ and $\X$) characterized by
$$P(\t,\X)\in \kK_{i,j}\iff \bdeg(P)=(i,j)\ \mbox{and} \ P(\t,\u(\t))=0.$$ 
It is known that 
(see for instance \cite{CSC98}) $$\kK_{*,1}=\oplus_{j=0}^\infty \kK_{j,1}\simeq\mbox{Syz}(\I)$$ is a free $\K[\t]$-module 
generated by two elements, one in degree $\mu$ for a positive integer $\mu$ such that $0\leq\mu\leq\frac{d}{2}$,
and the other of degree $d-\mu$. The identification with the syzygies is done via the obvious correspondence
\begin{equation}\label{corresp}
 \begin{array}{ccc}
\kK_{*,1}&\to&\mbox{Syz}(\I)\\
a(\t)X_1+b(\t)X_2+c(\t)X_3&\mapsto&\big(a(\t),b(\t),c(\t)\big).   
  \end{array}
\end{equation} 
For generic parametrizations of this type, we will have $\mu=\lfloor\frac{d}{2}\rfloor$, and it turns out that
$\mu=1$ if and only if the curve has a singularity of maximal rank (see for instance Corollary $1$ in \cite{CWL08}). 
We will focus here in the case $\mu=1$, note that this implies $d\geq2$. 
Assume w.l.o.g. - after a linear change of variables - that
\begin{equation}\label{K1} 
\kK_{*,1}=\langle t_1X_2-t_2X_1,a(\t)X_1+b(\t)X_2+c(\t)X_3\rangle_{\K[\t]}
\end{equation}
 with
$\deg(a)=\deg(b)=\deg(c)=d-1$. 
Let $\phi:\P^1\to\P^2$ be the map given by
\begin{equation}\label{param}
\begin{array}{ccc}
\P^1&\to&\P^2\\
(t_1:t_2)&\mapsto&(u_1(t_1,t_2):u_2(t_1,t_2):u_3(t_1,t_2))   
  \end{array}
\end{equation}
and set ${\mathrm C}:=\phi(\P^1)$. Note that $\phi$ is globally defined as we have assumed that the
$\gcd$ of the $u_i$'s is equal to $1$.
In the terminology of \cite{bus09}, $\kK$ is the moving curve ideal of the parametrization given by $\phi$. 
\par The following proposition holds straightforwardly, but will be very useful in the sequel. We define 
the \textit{degree} of the map $\phi$ given in (\ref{param}) as $\#\big(\phi^{-1}(\bf p)\big)$
for a generic point ${\bf p}\in{\mathrm C}$. The parametrization will be called \textit{proper} is its
degree is equal to one.
\begin{proposition}\label{birat}
If $t_1X_2-t_2X_1\in \kK_{1,1}$, then $\phi$ is a birational map and its inverse is given by
\begin{equation}\label{inverse}
\begin{array}{cccc}
\psi:&{\mathrm C}&\to&\P^1\\
&(X_1:X_2:X_3)&\mapsto&(X_1:X_2)   
  \end{array}
\end{equation}
\end{proposition}

\begin{proof}
As $u_3\neq 0$, then ${\mathrm C}$ is not included in the hyperplane $X_3=0$. Let us consider the affine
part of $\phi$ gotten by setting $t_2=1$ and $X_3=1$, so we have that ${\mathrm C}_{\mbox{aff}}:=\C\cap\{X_3=1\}$ is the image of
$$\begin{array}{cccc}
 \phi_{\mbox{aff}}:&{\mathcal U}&\to&\K^3\\  
&t_1&\mapsto&\left(\frac{u_1(t_1,1)}{u_3(t_1,1)},\frac{u_2(t_1,1)}{u_3(t_1,1)},1\right),
  \end{array}
$$where ${\mathrm U}:=\K\setminus\{t_1:\,u_3(t_1,1)=0\}$.
As $X_1-t_1X_2$ is a moving line that follows this parametrization, we have that 
$$\frac{u_1(t_1,1)}{u_3(t_1,1)}-t_1\frac{u_2(t_1,1)}{u_3(t_1,1)}=0,
$$ and from here we get that
$$
t_1=\frac{u_1(t_1,1)}{u_2(t_1,1)}=\frac{u_1}{u_2}
$$
in ${\mathrm V}:=\K\setminus\{u_2(t_1,1)u_3(t_1,1)=0\}$.
So we have that the following composition gives the identity in ${\mathrm V}$:
$$\begin{array}{cccccc}
 \phi_{\mbox{aff}}:&{\mathrm V}&\to&\K^3&\to&{\mathrm V}\\  
&t_1&\mapsto&\left(\frac{u_1(t_1,1)}{u_3(t_1,1)},\frac{u_2(t_1,1)}{u_3(t_1,1)},1\right)&&\\
&&&(X_1,X_2,X_3)&\to&\frac{X_1}{X_2}.
  \end{array}
$$
and hence $\phi_{\mbox{aff}}$ restricted to ${\mathrm V}$ has an inverse. In particular, this shows already
that $\deg(\phi)=1$.
\par
As explained in \cite[Exercise 7.8]{har92} (see also \cite{PSS02}), if $\mbox{char}(\K)=0$, then $\deg(\phi)=1$ implies
the birationality of $\phi$. In order to show the claim for any field of positive characteristic, due to 
\cite[page 77]{har92}, it is enough to show an open set ${\mathrm W}\subset \C$ such that
$${\mathrm W}\stackrel{\Psi}{\to}\P^1\stackrel{\phi}{\to}{\mathrm W}
$$
gives the identity. We choose $${\mathrm W}:=\{(x_1:x_2:x_3)\in\C:\,x_1x_2\neq0\}.$$ Note that
${\mathrm W}\neq\emptyset$,  otherwise we would have that every $\t\in\P^1$ is a zero of $u_1\,u_2$, which is impossible since the field $\K$ is infinite.
So we need to show that $$(x_1:x_2:x_3)=\phi\left(\psi(x_1:x_2:x_3)\right)$$ for all $(x_1:x_2:x_3)\in{\mathrm W}$, i.e.
$$(x_1:x_2:x_3)=\phi(x_1:x_2)=\left(u_1(x_1,x_2):u_2(x_1,x_2):u_3(x_1,x_2)\right).$$
This is equivalent to $x_iu_j(x_1,x_2)-x_ju_i(x_1,x_2)=0$ for all $i,j\in\{1,2,3\}$. As $\phi$ is defined globally,
there exists $\s:=(s_1,s_2)\in\K^2\setminus\{{\bf 0}\}$ such that $$(x_1:x_2:x_3)=(u_1(\s):u_2(\s):u_3(\s)).$$
So, we have 
$$x_i=\lambda_1^du_i(\s),\ i=1,2,3$$
for some $\lambda_1\in\K\setminus\{0\}$.
On the other hand, as $t_1X_2-t_2X_1\in \kK_{1,1}$, we also have
$$s_1x_2-s_2x_1=s_1u_2(s_1,s_2)-s_2u_1(s_1,s_2)=0.$$ So,
$$
(x_1:x_2)=(u_1(\s):u_2(\s))=(s_1:s_2),
$$
 and hence there exists $\lambda_2\in\K\setminus\{0\}$ such that
\begin{equation}\label{lig2}
\begin{array}{cclc}
x_i&=&\lambda_2\,s_i,&\ i=1,2,\\
u_j(x_1,x_2)&=&\lambda_2^d u_j(s_1,s_2),&\ j=1,2,3.
\end{array}
\end{equation}
From here, we get
$$x_iu_j(x_1,x_2)-x_ju_i(x_1,x_2)=\lambda_1^d\lambda_2^d(u_i(\s)u_j(\s)-u_j(\s)u_i(\s))=0,$$
which proves the claim.
\end{proof}

\begin{corollary}\label{cor}
 If a rational parametrization has $\mu=1$,  then the map $\phi$ is generically one-to-one onto $\C$.
\end{corollary}
As a consequence of this corollary, we can remove the  properness's requirement in the statements of both  Corollary $1$ 
in \cite{CWL08} and in
Theorem $2.3$ in \cite{CHW08}.

\begin{remark}
Curves (and more generally, hypersurfaces) having a singular point of maximal rank are known as \textit{monoid curves}
(resp. hypersurfaces) (see \cite{SZKD99,JLP08}). In our case, the fact that the map $\phi$ induces $\mu=1$ tell us not
only that ${\C}$ is a monoid curve but also that $\phi$ is a proper parametrization of it.
\end{remark}
\smallskip
\begin{definition}
A ``monoid parametrization'' of a rational curve is a proper rational map from $\P^1$ onto a monoid curve. Equivalently,
a monoid parametrization is a rational map $\phi:\P^1\to\P^2$ having $\mu=1$.
\end{definition}

All along this section we will work with a monoid parametrization $\phi$.
Let $E(\X)$ denote the irreducible equation of ${\mathcal C}$, which is defined modulo a nonzero constant in $\K$ and it is equal to 
$$E(\X)=\mbox{Resultant}_{\t}(t_1X_2-t_2X_1,a(\t)X_1+b(\t)X_2+c(\t)X_3).$$ 
As $\phi$ is proper, we know that this resultant gives the actual implicit equation and not a power of it. 
Computing it explicitely by using the \textit{Poisson Formula} (see for instance \cite{CLO05}), we get
\begin{lemma}\label{deg1}
Up to a nonzero constant, we have
 $$E(\X)=a(X_1,X_2)X_1+b(X_1,X_2)X_2+c(X_1,X_2)X_3.$$
In particular, $\deg_{X_3}(E(\X))=1$.
\end{lemma}

\begin{proposition}\label{prop1}
Let $F(\t,\X)\in\K[\t,\X]$. Then,
 $F(\t,\X)\in \kK$ if and only if $F(X_1,X_2,\X)$ is a polynomial multiple of $E(\X)$.
\end{proposition}
\begin{proof}
We can assume w.l.o.g. that $F(\t,\X)\in \kK_{i,j}$ for some $i,j\in\N$. 
Due to the irreducibility of $E(\X)$ and the Nullstellensatz, it is enough to show that if $(x_1:x_2:x_3)\in\C$
(i.e. $E(x_1,x_2,x_3)=0$),
then $F(x_1,x_2,x_1,x_2,x_3)=0$.
\par Let then $(x_1:x_2:x_3)\in\C$ and $(s_1:s_2)\in\P^1$ such that $(x_1:x_2:x_3)=\phi(s_1:s_2)$ (as $\phi$ is globally defined, we know that
$\mbox{im}(\phi)=\C$). So, there exists $\lambda\in\K\setminus\{0\}$
such that 
\begin{equation}\label{lig}
x_i=\lambda^du_i(s_1,s_2),\ i=1,2,3.
\end{equation}
Recall that  $F(\t,\X)\in \kK\iff F(\t,\u(\t))=0$, so we have
$$F(s_1,s_2,u_1(s_1,s_2),u_2(s_1,s_2),u_3(s_1,s_2))=0.$$
By using (\ref{lig}) and (\ref{lig2}) which also applies in this case, due to the fact that $F$ is 
homogeneous of bidegree $(i,j)$, we
then have that
$$\begin{array}{l}
F(x_1,x_2,x_1,x_2,x_3)=\lambda^{dj}F(x_1,x_2,u_1(s_1,s_2),u_2(s_1,s_2),u_3(s_1,s_2))
\\=\lambda^{dj}\lambda_2^iF(s_1,s_2,u_1(s_1,s_2),u_2(s_1,s_2),u_3(s_1,s_2))=0,
\end{array}
$$
and hence the claim holds.
\end{proof}

\begin{lemma}\label{lem1}
Let $F(\t,\X)$ be a bihomogeneous polynomial of bidegree $(i,j)\in\N^2$. Then there exists $G(\t,\X)$ bihomogeneous 
of bidegree $(i-1,i+j-1)$ such that
$$
X_2^iF(\t,\X)-t_2^iF(X_1,X_2,\X)=(t_1X_2-t_2X_1)G(\t,\X).
$$
\end{lemma}
\begin{proof}
$$
X_2^iF(\t,\X)-t_2^iF(X_1,X_2,\X)=F(t_1X_2,t_2X_2,\X)-F(t_2X_1,t_2X_2,\X).
$$
By applying the first order Taylor formula over the polynomial $p(\theta):=F(\theta,t_2X_2,\X)$, the claim follows
straightforwardly.
\end{proof}

We write as before $\kK=\oplus_{i,j} \kK_{i,j}$, $\kK_{i,j}$ being the space of 
\textit{moving curves of bidegree $(i,j)$ which
follow the parametrization (\ref{param})} (see \cite{CHW08} for more on this terminology).

\begin{proposition}\label{prop2}
For $(i,j)$ such that $i+j<d$, every nonzero element of $\kK_{i,j}$ is a polynomial multiple of $t_1X_2-t_2X_1$.
\end{proposition}

\begin{proof}
Let  $F(\t,\X)\in \kK_{i,j}$.
We have then
$$X_2^iF(\t,\X)-t_2^iF(X_1,X_2,\X)=(t_1X_2-t_2X_1)G(\t,\X)
$$
with -due to Proposition \ref{prop1}- $F(X_1,X_2,\X)$ a homogeneous polynomial multiple of $E(\X)$ of total
degree $i+j<d=\deg(E(\X))$. As $E(\X)$ is irreducible, we have then  $F(X_1,X_2,\X)=0$ and so
$$X_2^iF(\t,\X)=(t_1X_2-t_2X_1)G(\t,\X).
$$
The fact that $X_2$ is coprime with $t_1X_2-t_2X_1$ implies that the latter divides $F(\t,\X)$ and hence there exists
$G_0(\t,\X)$ such that
$F(\t,\X)=(t_1X_2-t_2X_1)G_0(\t,\X).
$
\end{proof}

\medskip
This result is optimal in the sense that we know that there is an element in $\kK$ of
degree $(d-1,1)$  which is not a multiple of $t_1X_2-t_2X_1$, namely 
$$a(\t)X_1+b(\t)X_2+c(\t)X_3\in \kK_{d-1,1}.$$ 
Also, we have $E(\X)\in \kK_{0,d}$ which clearly is not a multiple of $t_1X_2-t_2X_1$.
\par

\subsection{Construction of non trivial generators at degree $d-1$}

Now we will define one nonzero element in $\kK_{j,d-j}$ for $j=0,1,\ldots, d-1$. We will do this recursively starting 
from
$\kK_{d-1,1}$ and increasing the $\X$-degree. Set then $$F_{d-1}:=a(\t)X_1+b(\t)X_2+c(\t)X_3$$
and, for  $j$ from $1$ to $d-1$ do:
\begin{itemize}
\item write $F_{d-j}$ as $A_{d-j}(\t,\X)t_1+B_{d-j}(\t,\X)t_2$ (clearly there is more than one way of doing this, just choose one), 
\item Set $F_{d-j-1}:=A_{d-j}(\t,\X)X_1+B_{d-j}(\t,\X)X_2$.
\end{itemize}
We easily check that $F_j\in \kK_{j,d-j}$ for $j=0,\ldots, d-1$.
Note that this construction actually is the same as the one made in \cite{CHW08,bus09} with Sylvester forms. The simplicity of the 
moving line of degree one actually makes the whole determinantal presentation avoidable.

\begin{theorem}\label{mt}
$F_0,\ldots, F_{d-1}$ is a sequence of nonzero elements of $\K[\t,\X]$, moreover
$F_j(X_1,X_2,\X)=E(\X)$ for all $j=0,\ldots, d-1$. Let $\mathcal{J}$ be the ideal
 generated by ${\tt J}:=\{t_1X_2-t_2X_1, F_0,F_1,\ldots, F_{d-1}\}$. Then $\tt{J}$ is a minimal set of generators of 
$\mathcal{J}$.
\end{theorem}

\begin{proof}
From the construction of the $F_j$'s it is very easy to show inductively that $$F_j(X_1,X_2,\X)=E(\X).$$ This implies
in particular that none of the $F_j(\t,\X)$'s is zero.
\par Suppose now that one of them is a polynomial combination of the others. If $d=2$, we have 
 ${\tt J}=\{t_1X_2-t_2X_1,F_0,F_1\}$,
with $F_0$ of bidegree $(0,2)$ and the other two elements of bidegree $(1,1)$ being linearly  over
$\K[\t]$, as they are a basis of $\mbox{Syz}(\I)$. So, neither of them can be a polynomial combination of the others.
\par In the case, $d>2$, we have that
$t_1X_2-t_2X_1$ has total degree $2$ which is less than the
total degree of the $F_j$'s 
for all $j=0,\ldots,d-1$, so it cannot be a polynomial combination of the others. 
\par Suppose now
$$F_j(\t,\X)=G(\t,\X)(t_1X_2-t_2X_1)+\sum_{i\neq j}\lambda_i F_i(\t,\X),$$
with $\lambda_i\in\K$ (this is due to the fact that the total degree of all the $F_i$'s is $d$).
As $(i,d-i)\neq (j,d-j)$ if $i\neq j$, then we must have  $F_j(\t,\X)=(t_1X_2-t_2X_1)G_0(\t,\X)$, with
$G_0$ being the piece of bidegree $(i-1,j-1)$ of $G$. But this implies that
$$F_j(X_1,X_2,\X)=E(\X)\neq0.$$
\end{proof}

The following is the main result of this section.
\begin{theorem}\cite[Theorem $2.3$]{CHW08}\cite[Proposition $3.1$]{bus09}\label{maint1}
A minimal set of generators of $\kK$ is given by  $t_1X_2-t_2X_1, F_0,F_2,\ldots, F_{d-1}$.
\end{theorem}

\begin{proof}
 Following the notation of Theorem \ref{mt}, we have ${\mathcal J}\subset \kK$ and  the set 
$$\{t_1X_2-t_2X_1, F_1,F_2,\ldots, F_d\}$$
being a minimal family of generators of ${\mathcal J}$. The proof will follow if we show that ${\mathcal J}=\kK$.
Let $F(\t,\X)\in \kK_{i,j}$. If $i+j<d$, then thanks to Proposition \ref{prop2}, we know that $F(\t,\X)$
is a multiple of $t_1X_2-t_2X_1$, and hence an element of  ${\mathcal J}$.
Suppose now that $i+j\geq d$. By using Lemma \ref{lem1}, we get 
\begin{equation}\label{uno}
X_2^iF(\t,\X)-t_2^iF(X_1,X_2,\X)=(t_1X_2-t_2X_1)G(\t,\X)
\end{equation}
for some $G(\t,\X)\in\K[\t,\X]$. On the other hand, thanks to Proposition \ref{prop1},
 we know that $F(X_1,X_2,\X)$ is a polynomial multiple of $E(\X)$ i.e. there exists a polynomial $h(\X)$ of
degree $i+j-d$ such that $F(X_1,X_2,\X)=h(\X)E(\X)$. 

We will consider two different cases:
\begin{itemize}
 \item If $i<d$, we set $H(\t,\X):=h(\X)F_i(\t,\X)$. Clearly $H\in \kK_{i,j}$ and also,
$H(X_1,X_2,\X)=h(\X)E(\X)$. We also have, thanks to Lemma (\ref{lem1}),
\begin{equation}\label{dos}
 X_2^iH(\t,\X)-t_2^ih(\X)E(\X)=(t_1X_2-t_2X_1)\tilde{G}(\t,\X),
\end{equation}
for $\tilde{G}\in\K[\t,\X]$. We now substract (\ref{uno}) from (\ref{dos}) and get
$$X_2^i\big(F(\t,\X)-H(\t,\X)\big)=(t_1X_2-t_2X_1)G_0(\t,\X),
$$
with $G_0\in\K[\t,\X]$.
From here it is easy to see that $F(\t,\X)$ is a polynomial combination of $F_i$ and $t_1X_2-t_2X_1$.
\item If $i\geq d$, then note that as 
 $\deg_{X_3}\big(F(\t,\X)\big)=\deg_{X_3}\big(F(X_1,X_2,\X)\big)=j$ and,
as we saw  in Lemma \ref{deg1}, $\deg_{X_3}(E(\X))=1$. We then have that $\deg_{X_3}(h(\X))= j-1<i+j-d.$

On the other hand, as the total degree of $h(\X)$ is equal to $i+j-d$,
we then have that every monomial appearing in the expansion of $h(\X)$ has degree in the variables $X_1,X_2$ at least 
$i+j-d-j+1=i-d+1\geq1$.
We write then
$$h(\X)=\sum_{k=i-d+1}^{i-d+j}h_k(X_1,X_2)X_3^{i+j-d-k}$$
with $h_k(X_1,X_2)$ homogeneous of degree $k\geq i-d+1$.
For each $k$, let $H_k(\t,X_1,X_2)\in\K[\t,X_1,X_2]$ be any bihomogeneous polynomial of bidegree $(i-d+1,k-i+d-1)$ such that 
$H_k(X_1,X_2,X_1,X_2)=h_k(X_1,X_2)$. These are very easy to get, as we just  convert enough monomials in the $X_i$'s
in $h_k$ into $t_i$'s.
Now set
$$H(\t,\X):=\sum_{k=i-d+1}^{i-d+j}H_k(\t,X_1,X_2)X_3^{i+j-d-k}F_{d-1}(\t,\X).
$$
It is easy to see now that $H(\t,\X)\in \kK_{i,j}\cap {\mathcal J}$ and that $H(X_1,X_2,\X)=h(\X)E(\X)$. An argument like in
the previous case then shows that
$F(\t,\X)$ is equal to $H(\t,\X)$ plus a polynomial multiple of $t_1X_2-t_2X_1$, hence an element of ${\mathcal J}$.
\end{itemize}
\end{proof}

\smallskip
\begin{remark}\label{rem}
Note that if $i\geq d-1$ we can deduce from the proof of Theorem \ref{maint1} that actually 
$$\kK_{i,j}\subset \langle t_1X_2-t_2X_1,F_{d-1}(\t,\X)\rangle.$$ 
This fact is not surprising as -projectively- both the Rees Algebra
 and the Symmetric Algebra of $\I$ define the same scheme, and hence they are equal up to saturation. We 
will run into this situation again in the case of surfaces (proof of Theorem \ref{mmtt}).
\end{remark}

\bigskip
\section{Monoid hypersurfaces}\label{s2}
The elementary approach used in Section \ref{s1} can be extended straightforwardly to \textit{monoid hypersurfaces},
which are irreducible algebraic hypersurfaces having a singularity of multiplicity one less than their geometric degree 
(see \cite{SZKD99,JLP08}). Let $n\geq1$ be a positive integer, and consider then a monoid hypersurface
$\H\subset\P^{n}$ of degree $d$ having the point $(0:0:\ldots:0:1)\in\H$ of multiplicity $d-1$. It can be seen then
(see for instance \cite{JLP08}) that $\H$ can be rationally parameterized as follows:
$$\begin{array}{cccl}
\phi:&   \P^ {n-1}&\dasharrow&\P^{n}\\
&\ttheta&\mapsto&\left(\theta_1f_{d-1}(\ttheta):\ldots:\theta_{n}f_{d-1}(\ttheta):f_d(\ttheta)\right),
  \end{array}
$$
where $\ttheta:=(\theta_1:\ldots :\theta_{n})$ and $f_{d-1}(\t),\,f_d(\t)$ are homogeneous polynomials in
$\K[\t]=\K[t_1,\ldots,t_{n}]$ of respective degrees $d-1$ and $d$, without common factors. In this section, we will
set $\X=(X_1,\ldots, X_{n+1}).$ It is easy to 
see that the implicit equation of $\H$ is given by the polynomial 
\begin{equation}\label{EE}
E(\X)=f_{d-1}(X_1,\ldots,X_{n})X_{n+1}-f_d(X_1,\ldots,X_{n}).
\end{equation}
Also, it is straightforward to check that the inverse of $\phi$ is given by
$$\begin{array}{ccc}
  \H&\dasharrow&\P^{n}\\
(X_1:X_2:\ldots :X_{n+1})&\mapsto&(X_1:X_2:\ldots :X_{n})
  \end{array}
$$
As before, we can
then study the ideal of moving surfaces that follow this parameterization, i.e. the kernel of the map
$$\begin{array}{ccccc}
 \mathfrak{h}:&\K[\t,\X]&\to& \mbox{Rees}(\I)&\\
&X_i&\mapsto &t_if_{d-1}(\t)\,Z,&\ i=1,\ldots, n\\
&X_{n+1}&\mapsto& f_d(\t)\,Z,  
  \end{array}
$$
where $\I$ now stands for $\langle t_1f_{d-1}(\t),\,\ldots,t_{n}f_{d-1}(\t),\,f_d(\t)\rangle$, and
$\mbox{Rees}(\I):=\K[\t][\I\,Z]$.
Let $\kK$ be the kernel of $\mathfrak{h}$, which again happens to be a bigraded ideal, $\kK=\oplus_{i,j} \kK_{i,j}$.
We can distinguish the following elements in $\kK$:
\begin{itemize}
 \item $p_{i,j}(\t,\X):=t_iX_j-t_jX_i\in K_{1,1}$ for all $1\leq i< j\leq n$.
\item Write $f_d(\t)=t_1f_{d,1}(\t)+\ldots +t_{n}f_{d,n}(\t)$ and set
$$F_{d-1}(\t,\X):=f_{d-1}(\t)X_{n+1}-\sum_{j=1}^{n}f_{d,j}(\t)X_j.$$
Due to (\ref{EE}), we have that $F_{d-1}(X_1,\ldots,X_n,\X)=E(\X)$ and this shows that
$F_{d-1}(\t,\X)\in \kK_{d-1,1}\setminus\{\bf 0\}$.
\item Inductively, we write $F_j(\t,\X)\in \kK_{j,d-j}$ as $\sum_{i=1}^{n}F_{j,i}(\t,\X)t_i$, and then define
$$F_{j-1}(\t,\X):=\sum_{i=1}^{n}F_{j,i}(\t,\X)X_i.
$$
\end{itemize}
So, we have $F_{j-1}\in \kK_{j-1,d-j+1}\setminus\{\bf0\}$, having the property that
$$F_{j-1}(X_1,\ldots,X_n,\X)=E(\X).$$ This can be done for $j=d-1,d-2,\ldots,1,\ F_0(\t,\X)$ being equal to $E(\X)$.
\par
The following is the generalization of Theorem \ref{maint1} for monoid hypersurfaces.
\begin{theorem}\label{mt2}
 A minimal set of generators for $\kK$ is
$$\{p_{i,j},\,1\leq i<j\leq n,\,F_0,F_1,\ldots, F_ {d-1}\}.$$
\end{theorem}

\begin{proof}
As in Proposition \ref{prop1}, it is easy to see that $F(\t,\X)\in \kK$ if and only if $F(X_1,\ldots,X_{n},\X)$
is a multiple of the implicit equation. On the other hand, we again have that, if $F(\t,\X)$ has bidegree $(i,j)$,
then
$$X_{n}^iF(\t,\X)-t_{n}^iF(X_1,\ldots,X_n,\X)\in {\mathcal L}:=\langle p_{i,j}(\t,\X)\rangle.
$$
The key point here is that ${\mathcal L}$ is a prime ideal, as it is the ideal of $2\times 2$ minors
of the matrix
$$\left(\begin{array}{cccc}
t_1&t_2&\ldots &t_{n}\\
X_1&X_2&\ldots &X_{n}         
        \end{array}
\right),
$$
which is known to be prime (see for instance \cite{shar64}). So, as in the monoid curve case, 
it is easy to show that if $F(\t,\X)\in \kK_{i,j}$ with $i+j<d$, then $F(\t,\X)\in {\mathcal L}$. 
The case $i+j\geq d$ also
follows straightforwardly by using the same tricks as in the proof of Theorem \ref{maint1}. We leave the details to
the reader.
\end{proof}

\begin{remark}
Note that also in this case we have following property: $i\geq d-1$, then $\kK_{i,j}$ is again generated 
as a $\K[\t,\X]$-module by
the syzygies of $\I$ via a correspondence similar to (\ref{corresp})  (compare with Remark \ref{rem}).
\end{remark}

\bigskip
\section{Local Complete Intersection Surfaces with $\mu_1=\mu_2=1$}\label{s3}
Now we  focus in the case of some special parametric surfaces, so in this section se set $n=3$, $\t$ will stand for
$(t_1,t_2,t_3),\ \X$ will be $(X_1,X_2,X_3,X_4)$ and $\u=(u_1,u_2,u_3,u_4)$ a sequence of nonzero homogeneous 
polynomials of degree $d$ in $\K[\t]$. As before, let $\I:=\langle u_1,u_2,u_3,u_4\rangle\subset\K[\t]$. We also assume 
$\gcd(u_i(\t),\,i=1,2,3,4)=1$. We now look for generators of the kernel of the map
$$\begin{array}{clcc}
   \mathfrak{h}:&\K[\t,\X]&\to&\mbox{Rees}(\I)\\
&X_i&\mapsto&u_i\,Z
  \end{array}
$$
for $i=1,2,3,4$. As usual, $\mbox{Rees}(\I)=\K[\t][\I\,Z].$ Let $\kK\subset \K[\t,\X]$ be the kernel of $\mathfrak{h}$.

In this context, we do not have a nice $\mu$-basis situation like in the curve case. For instance, if $\K=\CC$ and
$V(\I)=\emptyset$, it is shown in \cite[Proposition 5.1]{cox01} that $\mbox{Syz}(\I)$ is not a free $\K[\t]$-module
anymore. 
So, in order
to be able to talk about $\mu$-bases again, we must require $V_{\P^2}(\I)\neq\emptyset$ plus one of 
the following conditions (which are
all equivalent to the fact that $\mbox{Syz}(\I)$ is a free graded module for any algebraically closed
field $\K$ , see \cite[Proposition 5.2]{cox01} for the case $\K=\CC$ and \cite{BCJ09} for a general $\K$:
\begin{itemize}
 \item the projective dimension of $\K[\t]/\I$ is $2$;
\item $\K[\t]/\I$ is Cohen-Macaulay;
\item $\langle t_1,t_2,t_3\rangle\notin Ass(\K[\t]/\I)$;
\item $\I$ is saturated with respect to the maximal ideal $\langle t_1,t_2,t_3\rangle$. 
\end{itemize}

If one (and all) of the above holds, we have an exact sequence of the form
\begin{equation}\label{res}
0\to \K[\t](-d-\mu_1)\oplus \K[\t](-d-\mu_2)\oplus \K[\t](-d-\mu_3)\to \K[\t](-d)^4\to \K[\t] \to \K[\t] / \I \to 0,
\end{equation}
where $\mu_i\in\N,\,\mu_1+\mu_2+\mu_3=d.$
\par
We will then restrict ourselves to this situation. We will consider in this section the case $\mu_1=\mu_2=1$. This forces $\mu_3=d-2$, so
we assume that $d>2$. Let then
$p_1(\t,\X),\,p_2(\t,\X),\,p_3(\t,\X)$ be a $\K[\t]$-basis of the Syzygy module of $\I$ which we regard again as linear forms
in the variables $\X$'s 
via the usual identification
$$\begin{array}{ccc}
\kK_{*,1}&\to&\mbox{Syz}(\I)\\
a(\t)X_1+b(\t)X_2+c(\t)X_3+d(\t)X_4&\mapsto&\big(a(\t),b(\t),c(\t),d(\t)\big),
\end{array}
$$
with $\deg_\t(p_i)=\mu_i$.
Write $$p_i(\t,\X)= L_{i1}(\X)t_1+L_{i2}(\X)t_2+L_{i3}(\X)t_3,\,i=1,2,
$$
with $L_{ij}(\X)$ being a homogeneous linear form in $\X$ with coefficients in $\K$. Set also
$M_j$ to be the signed $j$-th maximal minor of the matrix $\big(L_{ik}(\X)\big)_{1\leq i\leq 2, 1\leq k\leq 3}$, for $j=1, 2, 3$.

Consider now $\phi:\P^2\dasharrow\P^3$ given by
\begin{equation}\label{params}
\begin{array}{ccc}
\P^2&\dasharrow&\P^3\\
\t&\mapsto&\u(\t).   
  \end{array}
\end{equation}
Let $\S$ be the Zariski closure of $\phi(\P^1)$. 
\par The following proposition is the generalization of Proposition \ref{birat} for surfaces.
\begin{proposition}\label{birat2}
If $\I$ is saturated, $\mu_1=\mu_2=1$ and $\deg(\S)\geq3$, then $\phi$ has a rational inverse given by
\begin{equation}\label{iinverse}
\begin{array}{ccc}
{\S}&\dasharrow&\P^2\\
\X&\mapsto&(M_1:M_2:M_3).   
  \end{array}
\end{equation}
\end{proposition}

\begin{proof}
First, note that we cannot have all the $M_i$'s equal to zero, as this would imply that there exists $B(\X)$ and $A(\X)$
both nonzero such that
$$A(\X)p_1(\t,\X)=B(\X)p_2(\t,\X).
$$
But this is impossible as $\{p_1(\t,\X),p_2(\t,\X)\}$ is linearly  over $\K[\t]$, and the equality above 
would imply that they differ by a scalar in $\K$.
\par
It is also easy to see that $M_i(\X)t_j-M_j(\X)t_i\in \kK_{1,2}$ for all $i,j=1,2,3$. If $M_i(\X)=0$ for some $i,$ then
we will have $M_j(\u)u_i=0$ for $j\neq i$, i.e. $M_j(\u)=0$ for $j\neq i$. As one of these $M_j(\X)$'s must be 
different from zero, we would have an equation of degree $2$ vanishing over $\S$ which has degree at least three, a contradiction. Hence,
we actually have $M_i(\X)\neq0$ for all $i=1,2,3$.

As before, from the relations $ \frac{t_j}{t_i}=\frac{M_j(\u)}{M_i(\u)}$ in $\P^2\setminus\big(V(I)\cup 
\{t_1t_2t_3M_1(\u)M_2(\u)M_3(\u)=0\}\big),$ we have then that the following algebraic morphisms are inverse of each other
$$\begin{array}{ccccccc}
\P^2&\dasharrow&\S&&\S&\dasharrow&\P^2 \\
\t&\to&\u(\t)&&\X&\to&\underline{M},
\end{array}
$$
this completes the proof.
\end{proof}

\begin{example}\label{ex32}
Consider the parametrization given in \cite[Example $4.1$]{BCD03}, where the ideal $\I$ is saturated, having 
a Hilbert-Burch resolution of degrees $(\mu_1,\mu_2,\mu_3)=(1,1,2)$ with
$$\begin{array}{ccl}
p_1(\t,\X)&=&t_1X_1+t_1X_2+2t_2X_3+t_2X_4,\\
p_2(\t,\X)&=&2t_1X_1+t_2X_2+t_1X_3+3t_2X_4\\
p_3(\t,\X)&=&t_2t_3X_1+t_1^2X_2+t_2^2X_3+t_1t_3X_4.   
  \end{array}
$$
So, $u_1,\,u_2,\,u_3,\,u_4$ are the signed maximal minors of the following $4\times3$ matrix:
$$\left(\begin{array}{ccc}
t_1&2t_1&t_2t_3\\
t_1&t_2&t_1^2\\
2t_2&t_1&t_2^2\\
t_2&3t_2&t_1t_3
        \end{array}
\right).
$$
As $L_{13}(\X)=L_{23}(\X)=0$, we then have that $M_1(\X)=M_2(\X)=0$. The explanation here is that
$$M_3(\X)=\left|
\begin{array}{cc}
 X_1+X_2&2X_3+X_4\\ \\
2X_1+X_3&X_2+3X_4
\end{array}
\right|
$$ is the implicit equation of the parametric surface, which has degree $2$.
\end{example}

The following consequence is the generalization of Proposition \ref{prop1} for surfaces with $\mu_1=\mu_2=1$, 
its proof being straightforward.
\begin{corollary}\label{cor1}
If $\I$ is saturated, having a Hilbert-Burch resolution with degrees $(1,1,d-2)$ and $\deg(\S)\geq3$, then
$F(\t,\X)\in \kK$ if and only if $F(\M,\X)$ is a polynomial multiple of $E(\X)$ in $\K[\X]$.
\end{corollary}

\begin{corollary}\label{corr}
If in addition $V(\I)$ is a local complete intersection, then 
\begin{enumerate}
\item the implicit equation of $\S$ is equal to $p_3(\M,\X)$,
\item $\gcd(M_1,M_2,M_3)=1$.
\end{enumerate}

\end{corollary}
\begin{proof}
Let $E(\X)$ be the implicit equation of $\S$. We know that (see for instance \cite[Theorem $4.1$]{BCD03}) if $V(\I)$
is a l.c.i. then -up to a nonzero constant-
$$E(\X)=Res_{\tmu}(p_1,p_2,p_3).
$$ 
Here, $Res_\tmu$ stands for the \textit{homogeneous} resultant in the variables $\t$
associated to the sequence $\tmu=(1,1,d-2)$ as defined in \cite[Chapter 3]{CLO05}.
By using the Poisson formula for the homogenous resultant (see for instance \cite{CLO05}), it is easy to see in this situation 
that -again up to a nonzero constant-
$Res_{\tmu}(p_1,p_2,p_3)=p_3(\M,\X),$ the value of $p_3(\t,\X)$ in  $\big(M_1(\X):M_2(\X):M_3(\X)\big)$,
which is the only common zero of $p_1(\t,\X)$ and 
$p_2(\t,\X)$ in $\P^2_{\K(\X)}$. Hence, the first part of the statement follows. For the second item, if we
set $M_0:=\gcd(M_1,M_2,M_3)$ and $\deg(M_0)>0$, then we would have that
$$p_3(\M,\X)=M_0^{d-2}p_3\big(\frac{M_1}{M_0},\frac{M_2}{M_0},\frac{M_3}{M_0},\X\big)$$
with one of the two polynomials in the right hand side vanishing on $\S$. This is impossible as the degree of
each of them is strictly less than $$\deg(\S)=\deg\big(p_3(\M,\X)\big)=2d-3.$$
\end{proof}

\subsection{Minimal Generators of $\kK$}

We start with the following result which will be useful in the sequel.
\begin{proposition}\label{above}
Let ${\mathcal N}$ be the ideal generated by $p_1(\t,\X),p_2(\t,\X)$ in $\K[\t,\X]$. If $\deg(\S)\geq3,$ 
then ${\mathcal N}$ is a prime ideal.
\end{proposition}

\begin{proof}
Let ${\mathcal E}\subset\K[\X]$ be the ideal generated by $M_1(\X),\,M_2(\X),\,M_3(\X)$. Consider the following complex of $\K[\X]$ modules:
\begin{equation}\label{complex1}
0\to\K[\X]^ 2\stackrel{\psi_2}{\to}\K[\X]^3\stackrel{\psi_1}{\to}{\mathcal E}\to0
\end{equation}
with matrices
$$\begin{array}{ccl}
{\tt M}(\psi_1)&:=&\left(\begin{array}{ccc}
M_1(\X)&
M_2(\X)&
M_3(\X)
\end{array}\right)   
\\ \\
{\tt M}(\psi_2)&:=&\left(\begin{array}{cc}
L_{11}(\X)&L_{21}(\X)
\\L_{12}(\X)& L_{22}(\X)\\
L_{13}(\X)&L_{23}(\X)
\end{array}
\right).
  \end{array}
$$
The fact that $\gcd(M_1,M_2,M_3)=1$ (see for instance Corollary \ref{corr} (2))
tells us that $V_{\P^3}(M_1,M_2,M_3)$ has dimension at most $1$, as otherwise there would be an irreducible hypersurface 
${\mathrm W}\subset\P^3$ contained in this variety, and this would imply that the irreducible polynomial in $\K[\X]$
defining ${\mathrm W}$ divides each of the $M_i$, which is impossible. So we have then that the affine dimension of
$V({\mathcal E})$ is less than or equal to $2$, and hence
\begin{equation}\label{dosh}
\mbox{depth}({\mathcal E})=\mbox{codim}({\mathcal E})\geq 4-2=2,
\end{equation}
as $\K[\X]$ is a Cohen-Macaulay ring. This implies the converse of the Hilbert-Burch Theorem (see for instance
\cite[Theorem $1.4.17$]{BH93}), and we have that the complex (\ref{complex1}) is exact.
\par
Let ${\mathcal E}_1\subset\K[\X]$ be the ideal generated by $L_{ij}(\X),\,i=1,2,\,j=1,2,3$. We will show that
$\mbox{depth}({\mathcal E}_1)=\mbox{codim}({\mathcal E}_1)\geq3$.
If $V_{\P^3}({\mathcal E}_1)=\emptyset,$ then it is clear that the only prime containing ${\mathcal E}_1$ 
is the maximal ideal
$\langle X_1,X_2,X_3,X_4\rangle$ which has codimension $4$.
\par On the other hand, if $\underline{x}:=(x_1:x_2:x_3:x_4)\in V_{\P^3}({\mathcal E}_1)$, then we claim that 
$V_{\P^3}({\mathcal E}_1)=\{\underline{x}\}$.
Indeed if there is another point $\underline{y}:=(y_1:y_2:y_3:y_4)$ in this variety, consider 
${\tt S}_{\underline{x},\underline{y}}$, the $2$-dimensional 
$\K$-vector subspace generated by $ \{(x_1,x_2,x_3,x_4),\,(y_1,y_2,y_3,y_4)\}$ in $\K^4$, and 
${\tt S}_{\underline{x},\underline{y}}^{\bot}$, its ortoghonal
in $\big(\K^4\big)^*$, the latter will identified with the space of linear forms in $X_1,X_2,X_3,X_4.$
\par 
The fact that $\underline{x},\underline{y}\in V_{\P^3}({\mathcal E}_1)$
implies that
$L_{ij}(\X)\in {{\tt S}_{p,q}}^{\bot},\,i=1,2,\,j=1,2,3.$ After a linear change of variables, we may assume that
$X_1,X_2$ is a set of generators of ${\tt S}_{p,q}^{\bot}.$ As $M_1(\X),M_2(\X),M_3(\X)$ are quadrics in the 
$L_{ij}(\X)$'s with 
coefficients in $\K$ (recall that they are the signed maximal minors of  ${\tt M}(\psi_2)$), this implies that
the inverse of $\phi$ given in (\ref{iinverse}) in terms of the map $\M$ is actually a homogeneous polynomial 
function of two variables $X_1,\,X_2$
which covers a dense subset of $\P^2$, a contradiction. Hence, $V_{\P^3}({\mathcal E}_1)=\{\underline{x}\}$, 
and this means that 
$\mbox{codim}({\mathcal E}_1)=3$. So, we have in all the cases
\begin{equation}\label{unoh}
\mbox{depth}({\mathcal E}_1)\geq3.
\end{equation}

Note also that $p_1(\t,\X),\,p_2(\t,\X)$ is a regular sequence in $\K[\t,\X]$, as they are both irreducible
and of bidegree $(1,1)$. So, the ideal ${\mathcal N}$ is
unmixed. This, plus (\ref{dosh}) and (\ref{unoh}) put ourselves in  the conditions of \cite[Proposition $3.3$ (i)]{SV81},
which asserts that these conditions are equivalent to the fact that ${\mathcal N}$ is a prime ideal.
\end{proof}

\begin{remark}
The hypothesis $\deg(\S)\geq3$ is necessary, as one can see in Example \ref{ex32} that $M_3(\X)\,t_1\in {\mathcal N}$ but neither
$M_3(\X)$ nor $t_1$ are elements of this ideal.
\end{remark}

\smallskip
The case $\mbox{codim}({\mathcal E}_1)=3$ can actually happen, as the following example shows.

\begin{example}
Consider the following set of moving planes:
$$\begin{array}{ccl}
p_1(\t,\X)&=&t_1X_1+t_2X_2+t_3X_3,\\ 
p_2(\t,\X)&=&-t_1X_2+2t_2X_3-t_3X_1,\\
p_3(\t,\X)&=&t_1t_3X_1+t_1t_2X_2+t_2t_3X_3+t_2^2X_4.  
  \end{array}
$$
With the aid of the computer software {\tt Macaulay 2} (\cite{mac}) we check that the ideal generated by the signed maximal
minors of
$$\left(
\begin{array}{ccc}
t_1&-t_3&t_1t_3\\
t_2&-t_1&t_1t_2\\
t_3&2t_2&t_2t_3\\
0&0&t_2^2 
\end{array}
\right)
$$
is saturated and a local complete intersection in $\K[\t]$. Note that $(0:0:0:1)\in V({\mathcal E}_1),$
and it is the only element in this variety, as we saw in the proof of Proposition \ref{above}.
A basis of the syzygy module of this ideal is given by
$$\{p_1(\t,\X),\,p_2(\t,\X),\,p_3(\t,\X)\}.$$ The polynomial defining the implicit equation has degree $5=2*4-3$ and it is given by
$$\begin{array}{ccl}
E(\X)&=&X_1^ 4X_4  + 2X_1^3X_3^2  - 2X_1^3X_2X_3  - 4X_1^2X_3^3  - 2X_1X_2^2X_3^2 \\ &&- 2 X_1^2X_2X_3^2  - 2X_1X_2X_3^3 
+ X_1X_2^3X_3 - X_2^3X_3^2  - X_1^2X_2^3  + X_2^2X_3^2X_4  \\ &&- X_1^3X_2^2 - 2X_1^2X_2X_3X_4  + 2X_2^2X_3^3  + X_1^2X_2^2X_3.
\end{array}
$$
\end{example}
\subsection{Construction of generators of $\kK$}
Going back to our problem of constructing a set of minimal generators of $\kK$, we proceed as in the previous sections:
we will add to the syzygies $p_1(\t,\X),\,p_2(\t,\X),\,p_3(\t,\X)$ some elements of $\kK$ defined by specializing some 
variables $\t$'s into $\M(\X)$'s.
\begin{definition}
Set $F_{d-2}(\t,\X):=p_3(\t,\X)\in \kK_{d-2,1}$.
For $j$ from $d-2$ to $1$ do:
\begin{itemize}
 \item Write $F_j(\t,\X)\in K_{j,1+2(d-2-j)}$ as $$F_j(\t,\X)=A_j(\t,\X)t_1+B_j(\t,\X)t_2+C_j(\t,\X)t_3.$$ 
Again, there are several
ways of doing this, just choose any.
\item Set $F_{j-1}(\t,\X):=A_j(\t,\X)M_1(\X)+B_j(\t,\X)M_2(\X)+C_j(\t,\X)M_3(\X).$
\end{itemize}
Clearly, we have $F_j(\t,\X)\in \kK_{j,1+2(d-2-j)}$ for all $j=0,1,\ldots, d-2$. Also, it is clear that 
$F_j\neq0$ for all $j$,
as they all specialize to $F_0(\X)=E(\X)$ when we do $\t\mapsto\M$.
\end{definition}

\begin{theorem}\label{mmtt}
Let $\I:=\langle u_1,u_2,u_3,u_4\rangle$ 
be a saturated local complete intersection ideal generated by four nonzero polynomials homogeneous of degree 
$d$ in $\K[\t]$ with $d\geq3$,
such that a $\mu$-basis de $\I$ is given by three forms $p_1(\t,\X),p_2(\t,\X),p_3(\t,\X)$ 
of respective degrees $1,1,d-2$.
Then, a minimal set of generators of $\kK$ is \newline $\{p_1,\,p_2,\,F_0,F_1,\ldots, F_{d-2}\}$.
\end{theorem}

\begin{proof}
It is clear that $\I$ saturated + l.c.i implies  $\deg(\S)\geq3$. Indeed, due
to the degree formula (\cite{cox01}):
$$2d-3=\deg(\phi)\deg(\S),
$$
we get that actually $\deg(\S)=\deg(E(\X))$ is an odd number. If $\deg(\S)=1$, this implies that  the
polynomials $u_1,u_2,u_3,u_4$ are $\K$-linearly dependent, and this would imply the existence of a syzygy of
degree zero which contradicts our hypothesis. Hence, $\deg(\S)\geq3$ holds.
\par
Let ${\mathcal M}$ be the ideal generated by $p_1,\,p_2,\,F_0,F_1,\ldots, F_{d-2}$ in $\K[\t,\X]$. Clearly we have
${\mathcal M}\subset \kK$. As in the proof of Theorem \ref{mt}, let us show first that this set of generators is minimal.
\par We easily check that neither $p_1(\t,\X)$ nor $p_2(\t,\X)$ can be expressed as a polynomial combination
 of the others, on the other hand, if we have
$$F_j(\t,\X)=A(\t,\X)p_1(\t,\X)+B(\t,\X)p_2(\t,\X)+\sum_{i\neq j}C_i(\t,\X) F_i,
$$ if we pick the piece of bidegree $(j,1+2(d-2-j))$ in this expression, we get
$$F_j(\t,\X)=A_0(\t,\X)p_1(\t,\X)+B_0(\t,\X)p_2(\t,\X),
$$
which would imply that $Res_{(1,1,j)}(p_1,p_2,F_j)=0$, a contradiction with the fact that
the latter is a nontrivial multiple of $E(\X)$. So, the family of generators is minimal.
\par
Let us pick now an element $F(\t,\X)\in \kK_{i,j}$ for some $i,j\in \N$. Again we have to consider three cases:
\begin{itemize}
\item $2i+j<2d-3$: due to the bihomogeneities of $F$, we have as before
$$M_3^iF(\t,\X)-t_3^{i}F(\M,\X)\in {\mathcal N}=\langle p_1(\t,\X),\,p_2(\t,\X)\rangle.
$$
As we have $\deg\big(F(\M,\X)\big)=2i-j<2d-3$, and due to the fact that $F(\M,\X)$ must be a multiple of $E(\X)$ 
(see Corollary \ref{cor1}) which is irreducible of degree $2d-3$, we then have that $F(\M,\X)=0$ and hence
$M_3^iF(\t,\X)\in {\mathcal N},$ which is a prime ideal in virtue of Proposition \ref{above}. As $M_3$ is not in this ideal (it
has degree $0$ in the variables $\t$'s, and ${\mathcal N}$ is generated in $\t$-degree one), we then have 
$F(\t,\X)\in {\mathcal N}$.
\item $2i+j\geq 2d-3,\, i\leq d-2$: In this case we have $F(\M,\X)=h(\X)E(\X)$ with $h(\X)\in\K[\X]$. Consider then
$h(\X)F_i(\t,\X)\in \kK_{i,j},$ which actually satisfies
$$M_3^i\left(F(\t,\X)-h(\X)F_i(\t,\X)\right)\in {\mathcal N}.
$$
As before, due to the fact that ${\mathcal N}$ is prime and $M_3\notin {\mathcal N}$, we have then that
$$F(\t,\X)\in {\mathcal N}+\langle F_i(\t,\X)\rangle\subset {\mathcal M}.
$$
\item $2i+j\geq 2d-3,\,i\geq d-1:$ this case is actually be the most complicated one, as we do not have a way of
generalizing the proof given in Theorem \ref{mt} to this situation. However, deep results in Commutative Algebra
relating the Rees Algebra and the Symmetric Algebra of the ideal $\I$ can be applied in the case of $\I$ saturated and
local complete intersection. What follows can be found in \cite{BJ03,BCJ09} and the references therein:
\begin{enumerate}
 \item As $\I$ is a local complete intersection ideal, then it is of linear type outside 
$\mathfrak{m}:=\langle t_1,\,t_2,\,t_3\rangle$
and hence
$$\kK=\langle p_1(\t,\X),\,p_2(\t,\X),\,p_3(\t,\X)\rangle :{\mathfrak m}^\infty.$$
\item $$\kK / \langle p_1,\,p_2,\,p_3\rangle =\left(
\langle p_1,\,p_2,\,p_3\rangle :{\mathfrak m}^\infty\right) / \langle p_1,p_2,p_3\rangle \\
= H^0_{\mathfrak m}
\left(\K[\t,\X]/\langle p_1,p_2,p_3\rangle\right).$$
\item $\K[\t,\X]/\langle p_1,p_2,p_3\rangle$ is naturally identified with $\mbox{Sym}_{\K[\t]}(\I)$, the Symmetric 
Algebra of $\I$. 
\item As $\I$ is a saturated ideal,
due to \cite[Theorem $4$]{BC05}, we have that for  $\nu\geq d-2:$
$$H^0_{\mathfrak m}
\left(\K[\t,\X]/\langle p_1,p_2,p_3\rangle\right)_\nu=0.$$
Here, the grading is taken with respect to the $\t$-variables.
\end{enumerate}
This last result says that the class of $F(\t,\X)$ in $ \left(\kK / \langle p_1,\,p_2,\,p_3\rangle\right)_i$ is equal
to zero if $i\geq d-2$. Hence, we conclude that in our situation, 
as $$\deg_\t\big(F(\t,\X)\big)=i\geq d-1,$$ then $F(\t,\X)\in \langle p_1,\,p_2,\,p_3\rangle$. This completes the proof. 
\end{itemize}
\end{proof}

\end{document}